\def\N{\mathbb N}
\def\Z{\mathbb Z}
\def\R{\mathbb R}
\def\Q{\mathbb Q}
\def\A{\mathcal A}

\def\pfz{\begin{proof}}
\def\pfk{\end{proof}}

\documentclass[runningheads,a4paper]{article}
\usepackage{amssymb,amsmath,amsthm}
\usepackage{url}

\newtheorem{lem}{Lemma}
\newtheorem{thm}[lem]{Theorem}
\newtheorem{prop}[lem]{Proposition}

\newtheorem{ex}[lem]{Example}

\begin{document}


\title{Purely periodic expansions in systems with negative base}
\author{Z. Mas\'akov\'a\footnote{e-mail: zuzana.masakova@fjfi.cvut.cz}, E. Pelantov\'a\\[1mm]
{\normalsize Department of Mathematics FNSPE, Czech Technical University in Prague}\\
{\normalsize Trojanova 13, 120 00 Praha 2, Czech Republic}}
\maketitle
\begin{abstract}
We study the question of pure periodicity of expansions in the negative base numeration system.
In analogy of Akiyama's result for positive Pisot unit base $\beta$, we find a sufficient condition
so that there exist an interval $J$ containing the origin such that the $(-\beta)$-expansion of every rational number from $J$ is purely periodic.
We focus on the case of quadratic bases and demonstrate the following difference between the negative and positive bases: It is known that the finiteness property (${\rm Fin}(\beta)=\Z[\beta]$) is not only sufficient, but also necessary in the case of positive quadratic and cubic bases. We show that ${\rm Fin}(-\beta)=\Z[\beta]$ is not necessary in the case of negative bases.

%
\end{abstract}
%

\section{Introduction}

%
%
%
%
%

In every numeration system, one is interested in the description of numbers with finite or periodic representations. In classical positional systems with integer base, the
set of numbers with eventually periodic expansion coincides with $\Q$. If the base $\beta$ is irrational, then the set of numbers with periodic $\beta$-expansion contains also non-rational numbers, but in general, not all $\Q$. Schmidt~\cite{schmidt} observed that if every element of $\Q$ should have periodic $\beta$-expansion, then the base $\beta$ is a Pisot number or a Salem number. He also proved a surprising fact that for  $\beta$ root of $x^2-mx-1$, $m\geq 1$, all rational numbers from $[0,1)$ have purely periodic expansion. Since then, many results about purely periodic $\beta$-expansions appeared.
A very clear and detailed overview can be found in Adamczewski et al.~\cite{AFSS} and Akiyama et al.~\cite{AkiyamaBaratBertheSiegel}.

Our aim is to study similar questions in a numeration system with negative base, as defined by Ito and Sadahiro in~\cite{ItoSadahiro}. In order to follow the analogies for $(-\beta)$-expansions, we briefly recall here the main results concerning $\beta$-expansions. Define $\gamma(\beta)$ be the greatest number in $[0,1]$ such that the $\beta$-expansion of every rational $x\in[0,\gamma(\beta))$ is purely periodic. By an easy argument one can show that a necessary condition so that  $\gamma(\beta)>0$ is that $\beta$ is an algebraic integer. The result of Schmidt~\cite{schmidt} implies that such $\beta$ must be a Pisot or a Salem number. Akiyama~\cite{akiyamaGreedy} shows that $\beta$ must be an algebraic unit. He also gives an easy argument for the fact that $\beta$ cannot have positive conjugate, which, in turn, implies that $\beta$ is not a Salem number. Note also that Pisot numbers $\beta$ with positive conjugate cannot have the finiteness property, denoted by (F), namely that the set of finite $\beta$-expansions coincides with $\Z[\beta^{-1}]$, which for unit $\beta$ is equal to $\Z[\beta]$.
Altogether, $\gamma(\beta)>0$ implies that $\beta$ is a Pisot unit without positive conjugate.
As for sufficient condition, an important general result of Akiyama~\cite{akiyamaGreedy} states that
if $\beta$ is a Pisot unit satisfying (F) then $\gamma(\beta)>0$.

For the particular case of quadratic Pisot units one has, as a result of Schmidt~\cite{schmidt} and Hama and Iamahashi~\cite{HamaIamahashi}, that $\gamma(\beta)>0$ if and only if $\beta$ is a Pisot unit satisfying (F). Moreover, if this is the case, which happens if the conjugate of $\beta$ is negative, then $\gamma(\beta)=1$. In the opposite case with $\gamma(\beta)=0$ no rational in $[0,1]$ has purely periodic $\beta$-expansion.
Authors of~\cite{AFSS} study the case of cubic numbers and show that here
$\gamma(\beta)>0$ if and only if $\beta$ is a Pisot unit satisfying (F). Moreover, if
$\beta$ is not totally real, then the value $\gamma(\beta)$ is irrational.
In the case of $\beta$-expansions with quadratic and cubic $\beta$, the finiteness property (F) is therefore
not only sufficient, but also necessary for having $\gamma(\beta)>0$.

Let us now turn to the negative base systems, defined by a transformation $T_{-\beta}$ of the interval $I_\beta$. The first difference is that here 0 is contained in the interior of the interval $I_\beta$ and therefore one could ask about periodicity of expansions of positive and negative rationals in $I_\beta$ separately. Actually, such differentiation makes sense, as is clear from our results in the case when the base is quadratic. A reasonable analogue of the condition $\gamma(\beta)>0$ used for the positive case is now that there is a (not necessarily symmetric) interval $J\subset I_\beta$ containing 0 such that any rational $x\in J$ has purely periodic $(-\beta)$-expansion.

In this paper we first show -- following the arguments for positive base -- that a necessary condition for
having such an interval $J$ is that $\beta$ is a Pisot unit or Salem number (Section~\ref{sec:mainidea}).
We then show the analogue of Akiyama's sufficient condition for the existence of an interval $J$, namely that $\beta$ be a Pisot unit with finiteness property (Section~\ref{sec:general}) with infinite expansion of the left end-point of $I_\beta$.
  We focus on the case of quadratic Pisot units and describe the interval $J$ for both classes of them -- quadratic Pisot units with positive and negative conjugates. We show that for quadratic base $\beta$ an interval $J$ of pure periodicity exists, even if $\beta$ does not satisfy finiteness property (Sections~\ref{sec:kvadr1} and~\ref{sec:kvadr2}).

\section{Preliminaries}\label{sec:prelim}

Consider $\beta>1$. Following Ito and Sadahiro~\cite{ItoSadahiro}, we define $I_\beta=\big[-\frac{\beta}{\beta+1},\frac1{\beta+1}\big)$ and the $(-\beta)$-transformation $T_{-\beta}:I_\beta \to I_\beta$ by
$$
T_{-\beta} (x) = -\beta x - \Big\lfloor -\beta x + \frac{\beta}{\beta+1}\Big\rfloor\,.
$$
Using $T_{-\beta}$, one finds a representation of every $x\in I_\beta$ in the form
$$
x=\frac{x_1}{-\beta} + \frac{x_2}{(-\beta)^2} +\frac{x_3}{(-\beta)^3} +\cdots\,,\qquad x_i=\Big\lfloor -\beta T_{-\beta}^{i-1}(x) + \frac{\beta}{\beta+1}\Big\rfloor\in\A_\beta:=\{0,1,\dots,\lfloor\beta\rfloor\}\,,
$$
and we denote $d_{-\beta}(x):=x_1x_2x_3\cdots$. A string $y_1y_2y_3\cdots$ of positive integers is said to be {\em $(-\beta)$-admissible}, if there exists an $x\in I_\beta$ such that $d_{-\beta}(x)=y_1y_2y_3\cdots$.
A result of~\cite{ItoSadahiro} is that a string $y_1y_2y_3\cdots$ of positive integers is $(-\beta)$-admissible if and only if every suffix of it satisfies the following condition
\begin{equation}\label{eq:admis}
d_{-\beta}\big(\tfrac{-\beta}{\beta+1}\big) \preceq_{\text{alt}} y_i y_{i+1} y_{i+2}\cdots
\prec_{\text{alt}} \lim_{\varepsilon\to0+} d_{-\beta}\big(\tfrac{1}{\beta+1} - \varepsilon\big)\,,
\end{equation}
where the limit is taken over the product topology on $\A_{\beta}^{\N}$ and $\preceq_{\text{alt}}$ is the so-called alternate order on $\A_{\beta}^\N$. We write $x_1x_2x_3\cdots \prec_{\text{alt}} y_1y_2y_3\cdots$ if $(-1)^j(y_j-x_j)>0$ for the first index $j$ with $x_j\neq y_j$. Note also, that by~\cite{ItoSadahiro},
$$
\lim_{\varepsilon\to0+} d_{-\beta}\big(\tfrac{1}{\beta+1} - \varepsilon\big) = 0d_{-\beta}\big(\tfrac{-\beta}{\beta+1}\big)\,,
$$
unless $d_{-\beta}\big(\tfrac{-\beta}{\beta+1}\big)$ is purely periodic with odd period-length.

One can find a representation of every real number $x$ in the form $x= \sum_{i=1}^{\infty} x_i(-\beta)^{n-i+1}$ by first dividing $x$ by a suitable power of $(-\beta)$, so that
$\frac{x}{(-\beta)^{n+1}}\in I_{-\beta}$, finding $d_{-\beta}\big(\frac{x}{(-\beta)^{n+1}}\big)$
and multiplying back by $(-\beta)^{n+1}$. However, as explained for example in~\cite{ADMP},
such representation is not unique for countably many real numbers $x$. This phenomena is advantageously avoided when finding such power of $(-\beta)$ that $\frac{x}{(-\beta)^{n+1}}$ falls into the interior of
$I_\beta$. For, we can note that if $y\in\big(\frac{-\beta}{\beta+1},\frac1{\beta+1}\big)$, then $\frac{y}{-\beta}\in I_\beta$ and, moreover, $y+\frac{\beta}{\beta+1}\in(0,1)$, and therefore
$$
T_{-\beta}\Big(\frac{y}{-\beta}\Big) = -\beta \Big(\frac{y}{-\beta}\Big) - \left\lfloor-\beta \Big(\frac{y}{-\beta}\Big) +\frac{\beta}{\beta+1} \right\rfloor = y
\quad\text{and}\quad d_{-\beta}\Big(\frac{y}{-\beta}\Big) = 0d_{-\beta}(y)\,.
$$
Obviously, we also have
\begin{equation}\label{eq:transformident}
T_{-\beta}^j\Big(\frac{y}{(-\beta)^j}\Big) = y \quad\text{and}\quad d_{-\beta}\Big(\frac{y}{(-\beta)^j}\Big) = 0^jd_{-\beta}(y)
\quad\text{for all $j\geq1$.}
\end{equation}

It can be easily realized that every real $x$ can therefore be represented uniquely in the form
$x=\sum_{i=1}^{\infty} x_i(-\beta)^{n-i+1}$, where $x_1\neq0$ and the string $0x_1x_2x_3\cdots$
is $(-\beta)$-admissible. In agreement with~\cite{ADMP} we call this representation of $x$ the {\em $(-\beta)$-expansion} of $x$.

When the string $x_1x_2x_3\cdots$ ends in repeating zeros, denoted by $0^\omega$, we say that the $(-\beta)$-expansion of $x$ is finite. We denote the set of all real numbers with finite $(-\beta)$-expansion by ${\rm Fin}(-\beta)$. We also define the $(-\beta)$-integers as the numbers whose $(-\beta)$-expansion
does not contain negative powers of the base $(-\beta)$, i.e.
$$
\Z_{-\beta} = \left\{\sum_{i=1}^{n+1}x_i(-\beta)^{n-i+1} \,\Bigm|\, n\in\N_0,\, \text{ $0x_1x_2\cdots x_{n+1}0^\omega$ is a $(-\beta)$-admissible string}\right\}\,.
$$
Let us mention that the $(-\beta)$-integers were also studied in~\cite{steinerSubstRairo,dombek} from the combinatorial point of view; their geometric behaviour is in focus of~\cite{steiner}.

Obviously, every $x$ with finite $(-\beta)$-expansion can be written as a polynomial with integer coefficients in variables $\beta,\beta^{-1}$, formally ${\rm Fin}(-\beta)\subset\Z[\beta,\beta^{-1}]$.
In case that $\beta$ is an algebraic integer, this clearly reduces to
${\rm Fin}(-\beta)\subset\Z[\beta^{-1}]$. If, moreover, $\beta$ is an algebraic unit, one has
$\Z[\beta^{-1}]=\Z[\beta]$, and therefore ${\rm Fin}(-\beta)\subset\Z[\beta]$.

One is often interested whether the set of finite expansions is closed under arithmetic operations, i.e.
whether ${\rm Fin}(-\beta)=\Z[\beta,\beta^{-1}]$.
As shown in~\cite{MaPeVa}, it can happen only if $\beta$ is a Pisot number or a Salem number.
In the same paper, one identifies sofar the only known class of numbers with the above property, namely
 quadratic Pisot numbers with positive conjugate.


\section{Tools}\label{sec:mainidea}

Our aim is to describe $\beta$ for which there exists an interval $J\subset I_\beta$ such that all rational numbers from $J$ have purely periodic $(-\beta)$-expansion. The fact that a positive real
number $x$ has a purely periodic $(-\beta)$-expansion, i.e. $d_{-\beta}(x)=
(a_{N-1}\cdots a_0)^\omega$, means that
\begin{equation}\label{eq:purper}
x = z\Big( \frac1{(-\beta)^N} + \frac1{(-\beta)^{2N}} +\frac1{(-\beta)^{3N}}+\cdots \Big) =
\frac{z}{(-\beta)^N-1}\,,
\end{equation}
where $z=a_{N-1}(-\beta)^{N-1}+\cdots +a_1(-\beta) + a_0$. First, we shall
study necessary conditions under which such expression is possible for all
rational $x$ contained in an interval around the origin.

\begin{prop}\label{p:necessary}
Let $\beta>1$. Suppose there exists an interval $J$ containing zero such that the $(-\beta)$-expansion of any rational number in $J$ is purely periodic. Then $\beta$ is an algebraic integer, it is an algebraic unit, and it is either a Pisot or a Salem number. Moreover, if $\beta$ has a negative conjugate, then $J\subset(-\infty,0]$.
\end{prop}

\pfz
By assumption, there exists a $q\in\Z\setminus\{0\}$ such that $\frac1{q}\in J$, and thus it
has purely periodic $(-\beta)$-expansion. Without loss of generality, we may assume that the period length
is an even number $N\in\N$.
Necessarily,
 $\frac1{q} = z \frac1{\beta^N-1}$, where $z=\sum_{i=0}^{N-1}
 a_i(-\beta)^i$, $a_i\in\{0,1,\dots,\lfloor\beta\rfloor\}$, and thus
 $$
 \beta^N -1 = zq = \sum_{i=0}^{N-1} q a_i(-\beta)^i\,,
 $$
i.e. $\beta$ is a root of a monic polynomial with integer coefficients,
which implies that $\beta$ is an algebraic integer.

Suppose that $\beta$ is not a unit.
For sufficiently large $k$, we can put $q=\varepsilon N(\beta^k)$, where $N:\Q(\beta)\to\Q$ is the norm in the field $\Q(\beta)$ and $\varepsilon\in\{+1,-1\}$ is such that $\frac1{\varepsilon N(\beta^k)}\in J$. Thus the rational number
 $\frac1{\varepsilon N(\beta^k)}$ has purely periodic $(-\beta)$-expansion.
 Therefore
$\beta$ is a root of the polynomial
$$
\beta^N - \varepsilon N(\beta^k)\sum_{i=1}^{N-1}
 a_i(-\beta)^i - \varepsilon N(\beta^k)a_0-1 = 0\,.
$$
Obviously, $N(\beta)$ must divide the
constant term $- (\varepsilon N(\beta^k)a_0+1)$, which is not the case.

Now following the argument of Schmidt~\cite{schmidt} easily adapted to the negative
base in~\cite{MaPeISnumbers}, we can derive that $\beta$ cannot have conjugates of
modulus greater than $1$. In fact, one only
needs that there exists an interval containing $(-\beta)^{-j}$ whose rational numbers
have eventually periodic $(-\beta)$-expansion. This is by assumption true for sufficiently
large $j\in\N$, and thus we can conclude that $\beta$ is a Pisot or Salem number.

Let $x\in\Q\cap I_\beta$ have a purely periodic $(-\beta)$-expansion. Then by~\eqref{eq:purper}, we have
$x=\frac{z}{(-\beta)^N-1}$ for some $z=a_{N-1}(-\beta)^{N-1}+\cdots +a_1(-\beta) + a_0\in\Z_{-\beta}$,
where $N$ is the length of the period. Taking the field conjugate of $x$ in $\Q(\beta)$, we obtain
that
$$
x=x'=\frac{a_{N-1}(-\beta')^{N-1}+\cdots +a_1(-\beta') + a_0}{(-\beta')^N-1}\,.
$$
If the algebraic conjugate $\beta'$ of $\beta$ is negative, then the numerator of the above expression is
non-negative and the denominator is negative, thus necessarily $x\leq 0$.
\pfk

Note that the proof in fact implies that existence of a single positive rational $x$ with purely periodic $(-\beta)$-expansion implies that the conjugate of $\beta$ is not negative.

Let us turn our attention to a sufficient condition for $\beta$ guaranteeing the existence of the
desired interval $J$. Based on Proposition~\ref{p:necessary}, we may assume that $\beta$ is a Pisot
unit or a Salem number.
The sufficient condition must enable one to find to every $x\in J\cap \Q$ a $(-\beta)$-integer
$z$ satifying~\eqref{eq:purper} with, roughly speaking,
sufficiently ``short" $(-\beta)$-expansion, which concatenated infinitely many times gives a $(-\beta)$-admissible string. Let us summarize the points by which we shall proceed. They reflect that
$z=\big((-\beta)^N-1\big)x\in\Z_{-\beta}\subset {\rm Fin}(-\beta)\subset\Z[\beta]$.

\begin{itemize}

\item[(i)] For a rational number $x\in\Q(\beta)$ there exist
infinitely many $N\in\N$ such that
$z=\big((-\beta)^N-1\big)x\in\Z[\beta]$. This statement has been
mentioned in Akiyama~\cite{akiyamaGreedy} without proof. For reader's
convenience, we formulate and prove it as Lemma~\ref{l:apendix}.

\item[(ii)] For some base $-\beta$ and some $x$, we can derive that $z =
\big((-\beta)^N-1\big)x \in\Z[\beta]$ belongs to ${\rm
Fin}(-\beta)$. If $-\beta$ has finiteness property, i.e.
$\Z[\beta]={\rm Fin}(-\beta)$, this holds clearly for every
rational $x$. But it can happen also for $-\beta$ without
finiteness property, cf.~Lemma~\ref{l:finskladnoucarkou}.

\item[(iii)] For some base $-\beta$ and some $x$, we can derive that $z =
\big((-\beta)^N-1\big)x \in{\rm Fin}(-\beta)$ belongs to
$\Z_{-\beta}$, i.e.~there exists an expression for $z$ of the form
\begin{equation}\label{eq:betaintegr}
z= a_n(-\beta)^n + a_{n-1}(-\beta)^{n-1}+\cdots+a_1(-\beta)+a_0\,,
\end{equation}
where $0a_na_{n-1}\cdots a_1a_00^\omega$ is a
$(-\beta)$-admissible digit string and $a_n\neq0$. For that, we
use the estimate on the Galois image of $z$. This is done in
general in Lemma~\ref{l:obecodhad}. In the case of quadratic base,
we have exact description of the set $\Z_{-\beta}$ of
$(-\beta)$-integers, see Lemma~\ref{l:carky} and
Lemma~\ref{l:carky2}.

\item[(iv)]
If, moreover, $n<N$ and the string $(0^{N-1-n}a_na_{n-1}\cdots a_1a_0)^\omega$ is also $(-\beta)$-admissible, then
using
$\big((-\beta)^N-1\big)^{-1} = \sum_{j=1}^\infty(-\beta)^{-jN}$
we derive that the $(-\beta)$-expansion of $x$
is purely periodic, namely of the form
$$
d_{-\beta}(x)=d_{-\beta}\Big(\frac{z}{(-\beta)^N-1}\Big) =
(0^{N-1-n}a_na_{n-1}\cdots a_1a_0)^\omega\,,
$$
as we desired. For deciding about $n<N$, we use
Lemma~\ref{l:intervaly}. The question of admissibility is solved
in Lemma~\ref{l:obecadmis} for the general case. For quadratic
bases we derive stronger results using Lemmas~\ref{l:admisnetau}
and~\ref{l:admistau}.
\end{itemize}

Let us now present Lemmas~\ref{l:apendix} and~\ref{l:intervaly},
which will be recalled in solving both the general case and the
more precise results in the specific quadratic case.

\begin{lem}\label{l:apendix}
Let $\beta>1$ be an algebraic unit. Then for every $x\in\Q(\beta)$ there exists infinitely many even
$N\in\N$ such that $\big((-\beta)^N-1\big)x\in\Z[\beta]$.
\end{lem}

\pfz Every $x\in\Q(\beta)$ can be written in the form $x=\frac1q
z$, where $z\in\Z[\beta]$ and $q\in\N$. Since $\Z[\beta]$ is
closed under multiplication, it suffices to show the statement for
$z=\frac1q$. As $\beta$ is an algebraic integer, it is obvious
that $\beta^N\in\Z[\beta]$ for every $N\geq 0$, and from the ring
property of $\Z[\beta]$ we derive the implication
$$
(\beta^N-1)x\in\Z[\beta] \quad\implies\quad (\beta^{2N}-1)x=(\beta^N+1)(\beta^N-1)x \in\Z[\beta]\,.
$$
For every fix $q\in\N$ it is therefore sufficient to show
existence of a single $N$ with the property
$(\beta^N-1)\frac1q\in\Z[\beta]$.

Let $\beta$ be an algebraic unit of degree $d$ and let
$A\in\Z^{d\times d}$ be the companion matrix of $\beta$, i.e.
\begin{equation}\label{eq:matice}
A \left(\begin{smallmatrix} \beta^{d-1} \\[-2mm] \vdots \\[1mm] \beta \\ 1 \end{smallmatrix}\right) = \beta
\left(\begin{smallmatrix} \beta^{d-1}\\[-2mm] \vdots \\[1mm] \beta \\ 1 \end{smallmatrix}\right)
\end{equation}
On the set $\{A^k\mid k\in\N\}$ consider the equivalence relation
$A^k\sim A^l$ if $q$ divides all elements of the matrix $A^k-A^l$.
Obviously, the equivalence has finitely many classes, hence there
exists a class with at least two elements, i.e. indices $k>l$ such
that $A^k-A^l=qB$ for some integer matrix $B$. Therefore
$A^{k-l}-I=q(A^{-1})^lB$. Since $\beta$ is a unit, we have $\det A
= \pm 1$ and thus $A^{-1}\in\Z^{d\times d}$. Put $N=k-l$ and
$C=(A^{-1})^lB$. We have $N\in\N$, $C\in\Z^{d\times d}$ such that
\begin{equation}\label{eq:matice2}
A^N-I=qC\,.
\end{equation}
From~\eqref{eq:matice}, we derive
$$
qC \left(\begin{smallmatrix} \beta^{d-1} \\[-2mm] \vdots \\[1mm] \beta \\ 1 \end{smallmatrix}\right) =(\beta^N-1)
\left(\begin{smallmatrix} \beta^{d-1}\\[-2mm] \vdots \\[1mm] \beta \\ 1 \end{smallmatrix}\right)\,.
$$
Multiplying from the left by the row vector $(0,\dots,0,\frac1q)$, we obtain
$$
(0,\dots,0,1) \ C \left(\begin{smallmatrix} \beta^{d-1} \\[-2mm] \vdots \\[1mm] \beta \\ 1 \end{smallmatrix}\right) = (\beta^N-1)\frac1q\,.
$$
As $C\in\Z^{d\times d}$, the number on the left-hand side is an
integer combination of $1,\beta,\dots,\beta^{d-1}$, i.e.\ an
element of $\Z[\beta]$. This completes the proof. \
\pfk

The following Lemma~\ref{l:intervaly} allows one to estimate the
length of the $(-\beta)$-expansion of a $(-\beta)$-integer.

\begin{lem}\label{l:intervaly}
Let $z\in\R$. Denote by $n$ the most significant digit of $z$, i.e. such that
$$
z=x_1(-\beta)^n+x_2(-\beta)^{n-1}+x_3(-\beta)^{n-2}+\cdots\,,
$$
where $0x_1x_{2}x_{3}\cdots$ is a
$(-\beta)$-admissible digit string with $x_1\neq 0$.
$$
\begin{array}{lll}
\text{If} \quad \displaystyle{0<z<\ \frac{\beta^{2k+2}}{\beta+1}}\,,\ &\text{then $n$ is even and $n\leq 2k$.}\\[3mm]
\text{If} \quad \displaystyle{0>z>-\frac{\beta^{2k+1}}{\beta+1}}\,,\ &\text{then $n$ is odd and $n\leq 2k-1$.}
\end{array}
$$
\end{lem}

\pfz
From the assumption we have $d_{-\beta}(\frac{z}{(-\beta)^{n+1}})=x_1x_2x_3\cdots$ with
$0x_1x_2x_3\cdots$ admissible and therefore $\frac{z}{(-\beta)^{n+1}}\in\big(\frac{-\beta}{\beta+1},\frac1{\beta+1}\big)$. Since $x_1>0$,
we have $d_{-\beta}(\frac{z}{(-\beta)^{n+1}})=x_1x_2\cdots \prec_{\text{alt}} 0^\omega=d_{-\beta}(0)$. As the alternate
order on admissible strings corresponds to the natural order on reals, we derive that $\frac{z}{(-\beta)^{n+1}}<0$, and thus $n$ is even for $z>0$ and $n$ is odd for $z<0$.

Consider  $0<z< \frac{\beta^{2k+2}}{\beta+1}$. Then $n$ is even. Suppose $n\not\leq2k$, i.e. $n\geq 2k+2$. Set $y:=\frac{z}{(-\beta)^{2k+1}}$. We have
$$
\frac{-\beta}{\beta+1}<y=\frac{z}{(-\beta)^{2k+1}}<0\,,\quad\text{i.e.}\quad y\in\Big(\frac{-\beta}{\beta+1},\frac1{\beta+1}\Big)\,.
$$
Note that
$$
\frac{z}{(-\beta)^{n+1}} = \frac{y}{(-\beta)^j}\,,\quad \text{where }\ j=n+1-(2k+1)=n-2k\geq 2\,.
$$
Thus by~\eqref{eq:transformident}, we have
$$
d_{-\beta}\Big(\frac{z}{(-\beta)^{n+1}}\Big) = x_1x_2x_3\cdots = 0^jd_{-\beta}(y)\,,
$$
which is a contradiction with $x_1\neq 0$.

Let now $0>z> \frac{-\beta^{2k+1}}{\beta+1}$. Then $0<(-\beta)z<\frac{\beta^{2k+2}}{\beta+1}$
and $(-\beta)z = x_1(-\beta)^{n+1}+x_2(-\beta)^n+\cdots$ with $0x_1x_2x_3\cdots$ admissible and $x_1\neq0$.
By the first part of the proof, $n+1$ is even and $n+1\leq 2k$. This implies $n$ is odd and $n\leq 2k-1$.
\pfk

%
%

%
%
%

\section{Quadratic case $\boldsymbol{\beta^2=m\beta-1}$}\label{sec:kvadr1}

In this section we study the case when the base $-\beta$ satisfies $\beta^2=m\beta-1$. Crucial in deriving Theorem~\ref{thm6} is the knowledge of admissible digit strings in this numeration system. For that we use the admissibility condition~\eqref{eq:admis} with
$$
d_{-\beta}\big(\tfrac{-\beta}{\beta+1}\big)=\big((m-1)1\big)^\omega\qquad \lim_{\varepsilon\to 0+}
d_{-\beta}\big(\tfrac{1}{\beta+1}-\varepsilon\big) = 0\big((m-1)1\big)^\omega\,.
$$
The admissibility condition~\eqref{eq:admis} for the string $y_1y_2y_3\cdots$ now thus reads
$$
\big((m-1)1\big)^\omega \preceq_{\text{\tiny alt}} y_{i}y_{i+1}\cdots \prec_{\text{\tiny alt}} 0\big((m-1)1\big)^\omega\,.
$$
It is not difficult to see that a string $y_1y_2y_3\cdots$ of digits in $\{0,1,\dots,m-1\}$ is admissible if and only if
\begin{equation}\label{eq:admisquadrposit}
(i)\quad y_{i}=m-1 \ \Rightarrow \ y_{i+1}\geq 1\qquad \text{ and }\qquad (ii)\quad \text{it does not end with the suffix $0((m-1)1)^\omega$}.
\end{equation}
With this, one can easily verify that for example the $(-\beta)$-expansion of $-\tfrac12$ is
of the following form
$$
d_{-\beta}(-\tfrac12) =
\begin{cases}
 (kk1)^\omega & \text{if $m=2k-1$, $k\geq 2$,}\\
 (k1)^\omega & \text{if $m=2k$, $k\geq 2$.}\\
\end{cases}
$$
As an example of a positive rational number, one can take $\frac1{m+1}$, which satisfies
$\frac1{m+1}< \frac1{\beta+1}$ for all $m\geq 3$. We have
$$
d_{-\beta}(\tfrac1{m+1}) =\big(0(m-1)1110\big)^\omega\,.
$$

In fact, as shown in the following theorem, all rational numbers in the interval $I_{\beta}$ have purely periodic expansion.

\begin{thm}\label{thm6}
Let $\beta$ be a quadratic Pisot unit with positive conjugate. Then every rational $x\in I_\beta$ has purely periodic $(-\beta)$-expansion.
\end{thm}

The proof will use two auxiliary results (following the guideline in Section~\ref{sec:mainidea}).

\begin{lem}\label{l:carky}
Let $\beta>1$ satisfy $\beta^2=m\beta-1$, $m\geq 3$. Then
$$
\Z_{-\beta} = \big\{z\in\Z[\beta] \mid z' \in (-\tfrac1\beta H,
H)\big\}\,,
$$
where $H=\beta \frac{\beta-1}{\beta+1}$ and $z'$ denotes the image
of $z$ under the Galois automorphism of the field $\Q(\beta)$.
\end{lem}

\pfz We use constants
$$
H=\frac{1-\beta'}{\beta'(1+\beta')}=\beta
\frac{\beta-1}{\beta+1}\quad\text{and}\quad
K=\frac{\beta'(1-\beta')}{1+\beta'}=\frac{1}{\beta^2}H
$$
from the proof of Theorem~13 in~\cite{MaPeVa}, where it is shown
that
$$
-\frac1\beta H<z'<H \ \text{ for any
}z\in\Z_{-\beta},\qquad\text{and}\qquad z'>K \ \text{ for any }
z\in\Z_{-\beta}\setminus(-\beta)\Z_{-\beta}\,.
$$
The inclusion $\subset$ follows immediately.

In order to show that any $z\in\Z[\beta]$ satisfying
$z'\in(-\tfrac1\beta H, H)$ is a $(-\beta)$-integer, recall
from~\cite{MaPeVa} that for this base we have $\Z[\beta]={\rm
Fin}(-\beta)$. Suppose for contradiction that
$z=\sum_{i=s}^na_i(-\beta)^i$ with $a_{s}\neq 0$ for some $s<0$.
Thus $z=(-\beta)^{s}w$ for some
$w\in\Z_{-\beta}\setminus(-\beta)\Z_{-\beta}$. We have
$w'>K$. If $s$ is even, i.e. $-s\geq 2$, then
$$
z'=\big((-\beta)^{s}w\big)'=\beta^{-s}w'>\beta^{-s}K\geq
\beta^2K=H\,,
$$
which is a contradiction. If $s$ is odd, i.e. $-s\geq 1$, we have
$$
z'=\big((-\beta)^{s}w\big)'=-\beta^{-s}w'<-\beta^{-s}K\leq
-\beta K=-\tfrac1\beta H\,,
$$
which is again a contradiction.
\pfk

\begin{lem}\label{l:admisnetau}
Let $\beta>1$ satisfy $\beta^2=m\beta-1$, $m\geq 3$. Then
admissibility of the string $0x_1x_2\cdots x_k0^\omega$ implies
admissibility of $(0^tx_1x_2\cdots x_k)^\omega$ for any $t\geq 0$.
\end{lem}

\pfz
This is not difficult to see, realizing that the admissibility condition for the string $y_1y_2y_3\cdots$
is given by~\eqref{eq:admisquadrposit}.
We can derive that admissibility of the given string
$0x_1x_2\cdots x_k0^\omega$ implies that $x_k\leq m-2$. It is then
easy to see that also  $(0^tx_1x_2\cdots x_k)^\omega$ is
admissible. \pfk

\pfz[Proof of Theorem~\ref{thm6}]
Consider a non-zero $x\in\Q\cap I_\beta = \big[-\frac{\beta}{\beta+1},\frac1{\beta+1}\big)$. By Lemma~\ref{l:apendix} there exists infinitely many even $N$ such that $z=\big((-\beta)^N-1\big)x\in\Z[\beta]$. As $|-\beta'|<1$ and $x\neq-\frac{\beta}{\beta+1}\notin\Q$, one can choose $N$ sufficiently large so that
$$
z' = \big((-\beta')^N-1\big)x \in (-r,-l)=\big(-\tfrac1{\beta+1},\tfrac{\beta}{\beta+1}\big)\subset\big(-\tfrac1\beta H,H\big)\,,
$$
where $H=\beta \frac{\beta-1}{\beta+1}$. According to Lemma~\ref{l:carky}, we have $z\in\Z_{-\beta}$, and thus
it can be written in the form
\begin{equation}\label{eq:cislo4}
z = a_n(-\beta)^n + a_{n-1}(-\beta)^{n-1}+\cdots+a_1(-\beta)+a_0\,,
\end{equation}
where $0a_na_{n-1}\cdots a_1a_00^\omega$ is a $(-\beta)$-admissible digit string and $a_n\neq0$.
If $x>0$, we have
$$
0<z<\frac{(-\beta)^N-1}{\beta+1}<\frac{(-\beta)^N}{\beta+1}
$$
Since $N$ is even, by Lemma~\ref{l:intervaly}, the position $n$ of the first significant digit of $z$ satisfies $n\leq N-2<N-1$.
If $x<0$, we have
$$
0>z=\big((-\beta)^N-1\big)x>\frac{\big((-\beta)^N-1\big)(-\beta)}{\beta+1}>\frac{(-\beta)^{N+1}}{\beta+1}
$$
Since $N+1$ is odd, by Lemma~\ref{l:intervaly}, the position $n$ of the first significant digit of $z$ satisfies $n\leq N-1$.

Now by Lemma~\ref{l:admisnetau}, the string
$(0^{N-1-n}a_na_{n-1}\cdots a_1a_0)^\omega$ is admissible, and so
by item (iv) in Section~\ref{sec:mainidea}, the
$(-\beta)$-expansion of $x$ is purely periodic.
\pfk

\section{Quadratic case $\boldsymbol{\beta^2=m\beta+1}$}\label{sec:kvadr2}

With analogy to the positive base systems, where either every or no rational in $[0,1)$ has purely periodic expansions, one could expect that having purely periodic expansions of all rationals in $I_\beta$ for quadratic Pisot unit $\beta$ with positive conjugate, then if $\beta$ has negative conjugate, no rational in $I_\beta$ would have purely periodic expansion.

We can even support this idea by an example of a rational number in $I_\beta$, e.g. $\frac1{2m}$, whose expansion is given by
$$
d_{-\beta}(\tfrac1{2m})=\begin{cases}
d_{-\beta}(\tfrac1{4k-2}) = 0k\big((k-1)k(2k-1)(2k-1)(2k-2)(k-1))\big)^\omega  &\text{if $m=2k-1$, $k\geq 2$,}\\
d_{-\beta}(\tfrac1{4k}) = 0(k+1)\big((2k)(2k)(2k-1)k)\big)^\omega  &\text{if $m=2k$, $k\geq 1$},
\end{cases}
$$
Note that $\frac1{2m}$ satisfies $\frac1{2m}< \frac1{\beta+1}$ only $m\geq 2$.
For the case $m=1$, i.e. $\beta=\frac12(1+\sqrt5)$, the golden ratio, consider for example
$$
d_{-\beta}(\tfrac13) = 01(00111100)^\omega\,.
$$

In fact, it turns out that taking any positive rational in $I_\beta$, the corresponding $(-\beta)$-expansion
is indeed not purely periodic.
However, taking an example of a negative rational convinces us that some $x\in I_\beta\cap \Q$ may have purely periodic expansions. Consider $-\frac12$. We have
$$
d_{-\beta}(-\tfrac12) =
\begin{cases}
 (k(k-1)0)^\omega & \text{if $m=2k-1$, $k\geq 1$,}\\
 (k0)^\omega & \text{if $m=2k$, $k\geq 1$.}
\end{cases}
$$

In order to study this question, recall that we now have
$$
d_{-\beta}\big(\tfrac{-\beta}{\beta+1}\big)=m(m-1)^\omega\qquad \lim_{\varepsilon\to 0+}
d_{-\beta}\big(\tfrac{1}{\beta+1}-\varepsilon\big) = 0m(m-1)^\omega\,.
$$
The admissibility condition~\eqref{eq:admis} for the string $y_1y_2y_3\cdots$ now thus reads
$$
m(m-1)^\omega \preceq_{\text{\tiny alt}} y_{i}y_{i+1}\cdots \prec_{\text{\tiny alt}} 0m(m-1)^\omega
$$
and it is not difficult to see that a string $y_1y_2y_3\cdots$ of digits in $\{0,1,\dots,m\}$
is admissible if and only if
\begin{equation}\label{eq:admisquadrnegat}
\begin{array}{cl}
(i)&\text{it does not contain the
substring $m(m-1)^{2k}A$, where $A\leq m-2$, $k\geq 0$,}\\
(ii)&\text{it does not contain the
substring $m(m-1)^{2k+1}m$, $k\geq 0$,} \\
(iii)&\text{it does not end in $0m(m-1)^\omega$.}
\end{array}
\end{equation}

Using this, we show that if the base $-\beta$ satisfies $\beta^2-m\beta-1$, then positive and negative rational numbers have different behaviour with respect to $(-\beta)$-expansions.

\begin{thm}\label{thm:tau}
Let $\beta$ be a quadratic Pisot unit with negative conjugate. Then every rational $x\in I_\beta\cap(-\infty,0]$ has purely periodic $(-\beta)$-expansion
and no rational $x\in I_\beta\cap (0,+\infty)$ has purely periodic $(-\beta)$-expansion.
\end{thm}

We would like to write similar auxiliary results as in the previous section, namely Lemmas~\ref{l:carky} and~\ref{l:admisnetau}. The analogue of the latter concerning admissibility strings is formulated as Lemma~\ref{l:admistau}. However, the proof of Lemma~\ref{l:carky} uses the fact that for quadratic Pisot units with positive conjugates ${\rm Fin}(-\beta)$ is a ring and therefore ${\rm Fin}(-\beta)=\Z[\beta]$.
In case that $\beta$ is a root of $x^2-mx-1$, the set ${\rm Fin}(-\beta)$ is not a ring, since for example the $(-\beta)$-expansion of $-1$ is infinite, and thus necessarily, one has $\Z[\beta] \supsetneq {\rm Fin}(-\beta)$. Nevertheless, an analogue of Lemma~\ref{l:carky} can still be stated, see Lemma~\ref{l:carky2}.
For, as was shown in~\cite{MaVaKybernetika}, ${\rm Fin}(-\beta)$ is closed under addition. Consequently, every element $x$ of the form $x=a+b(-\beta)$ for $a,b\in\Z$, $a,b\geq 0$, has finite $(-\beta)$-expansion. In fact, one can characterize the elements of
$\Z[\beta]$ that also belong to ${\rm Fin}(-\beta)$ using their Galois conjugate.

\begin{lem}\label{l:finskladnoucarkou}
Let $\beta>1$ satisfy $\beta^2=m\beta+1$, $m\geq 1$. Then
$$
{\rm Fin(-\beta)}=\{z\in\Z[\beta]\mid z'\geq 0\}\,.
$$
\end{lem}

\pfz
Take $z=\sum_{i=s}^ka_i(-\beta)^i\in{\rm Fin}(-\beta)$. Since $\beta'=\frac{-1}{\beta}$, the Galois image of $z$ satisfies $z'=\sum_{i=s}^ka_i\beta^{-i}\geq 0$.

For the opposite inclusion, consider $z=a+b(-\beta)$, $a,b\in\Z$, such that $z'=a+b\beta^{-1}\geq 0$.
%
%
By induction, one easily sees that
\begin{equation}\label{eq:g}
\frac{z}{(-\beta)^k} = G_{k+1} + G_{k}(-\beta)\,,
\end{equation}
where the sequence $(G_k)_{k\geq 0}$ satisfies the recurrence
$$
G_0=b\,,\quad G_1=a\,,\qquad G_{k+2}=mG_{k+1}+G_k\,.
$$
Since $z'>0$, we have $\big(\frac{z}{(-\beta)^k}\big)'=z'\beta^k = G_{k+1} + G_{k}\beta^{-1}>0$, and thus two consecutive elements
of the sequence $(G_k)_{k\geq 0}$ are either both non-negative or have opposite signs. We show that it is not possible that $G_kG_{k+1}<0$ for all $k\geq 0$. Suppose, without loss of generality, that $G_{2j}>0$ and $G_{2j+1}<0$ for all $j\geq 0$.
From the recurrence, we have
$$
0<G_{2j}=mG_{2j-1}+G_{2j-2} < G_{2j-2}\,,
$$
i.e. $(G_{2j})_{j\geq0}$ is a strictly decreasing sequence of
positive integers, which is impossible. Thus, there exists an
index $k$ such that both $G_k$ and $G_{k+1}$ are non-negative.
Therefore $\frac{z}{(-\beta)^k} = G_{k+1} + G_{k}(-\beta)$ is a
sum of a finite number of elements $1,-\beta$ of ${\rm
Fin}(-\beta)$, which (by the result of~\cite{MaVaKybernetika})
implies that $z\in{\rm Fin}(-\beta)$. \pfk

\begin{lem}\label{l:carky2}
Let $\beta>1$ satisfy $\beta^2=m\beta+1$, $m\geq 1$. Then
$$
\Z_{-\beta} = \big\{z\in\Z[\beta] \mid z' \in [0, \beta)\big\}\,.
$$
where $z'$ denotes the image of $z$ under the Galois automorphism
of the field $\Q(\beta)$. In particular, for $z\in\Z[\beta]$ we
have $z\in(-\beta)\Z_{-\beta}$ if and only if $z'\in[0,1)$.
\end{lem}

\pfz Take a $(-\beta)$-integer $z=\sum_{i=0}^na_i(-\beta)^i$,
$a_i\in\{0,1,\dots,m\}$. Obviously $z\in{\rm Fin}(-\beta)
\subset\Z[\beta]$. Using $\beta'=-\tfrac1\beta$ we estimate
$$
0\leq z' = \sum_{i=0}^n\frac{a_i}{\beta^i} < m-1 +
\sum_{i=1}^{+\infty} \frac{m}{\beta^i} = \frac{m}{1-\beta^{-1}} -1
= \beta\,,
$$
where we have used that an admissible string $a_n\cdots
a_00^\omega$ cannot have $a_0=m$. The inclusion $\subset$ is thus
proved.

Suppose on the other hand that $z\in\Z[\beta]$ satisfies
$z'\in[0,\beta)$. Since $z'$ is non-negative, by
Lemma~\ref{l:finskladnoucarkou}, $z\in{\rm Fin}(-\beta)$ and
therefore $z$ is of the form $z=\sum_{i=s}^na_i(-\beta)^i$ with
$a_{s}\neq 0$ for some $s\in\Z$. Then
$$
\beta> z' =\sum_{i=s}^n\frac{a_i}{\beta^i} \geq \beta^{-s}\,,
$$
which implies $s\geq 0$, i.e. $z\in\Z_{-\beta}$.
\pfk

\begin{lem}\label{l:admistau}
Let $\beta>1$ satisfy $\beta^2=m\beta+1$, $m\geq 1$. If
$0x_1x_2\cdots x_k0^\omega$ is admissible, then  is
$(0^tx_1x_2\cdots x_k0)^\omega$ admissible for any $t\geq 0$.
\end{lem}

\pfz It is not difficult to see that if the string $0x_1x_2\cdots x_k0^\omega$
satisfies conditions (i)--(ii) from~\eqref{eq:admisquadrnegat}, then the same is true
for the string $(0^tx_1x_2\cdots x_k0)^\omega$, which is to be shown.
\pfk

\pfz[Proof of Theorem~\ref{thm:tau}]
By Proposition~\ref{p:necessary}, no positive $x\in\Q\cap I_\beta$ has a purely periodic $(-\beta)$-expansion.

Let us now show that every rational $x\in\big(-\frac{\beta}{\beta+1},0\big)$ has purely periodic $(-\beta)$-expansion. By Lemma~\ref{l:apendix}, there exists
a sufficiently large even $N$ such that $z=\big((-\beta)^N-1\big)x\in\Z[\beta]$ and
$$
z'=\big((-\beta')^N-1\big)x \in \Big(0,\frac{\beta}{\beta+1}\Big) \subset (0,1)\,.
$$
Therefore, by Lemma~\ref{l:carky2}, we may write $z= a_n(-\beta)^n + \cdots + a_1(-\beta) + a_0$ with $a_0=0$.
Since
$$
0>z>\frac{\big((-\beta)^N-1\big)(-\beta)}{\beta+1}>\frac{(-\beta)^{N+1}}{\beta+1}
$$
and $N+1$ is odd, by Lemma~\ref{l:intervaly}, the position $n$ of
the most significant digit of $z$ satisfies $n\leq N-1$. In order
to conclude pure periodicity of the $(-\beta)$-expansion of $x$,
realize that the string $(0^{N-n+1}a_n\cdots a_0)^\omega$ is
admissible by Lemma~\ref{l:admistau} with the use of $a_0=0$. \pfk

\section{General case}\label{sec:general}

Let us move to the study of a sufficient condition for existence of interval $J$ of pure periodicity
for general bases $-\beta$. By Proposition~\ref{p:necessary}, we may assume that
$\beta$ is a Pisot unit or a Salem number. The key is to put a stronger assumption, namely that $\beta$ satisfies the finiteness property, i.e. ${\rm Fin}(-\beta)=\Z[\beta]$. From~\cite{MaPeVa} we know that this happens only for Pisot or Salem numbers. We also know (cf.~\cite{ChiaraFrougny}) that if $\beta$ is a Pisot number, than the $(-\beta)$-expansion of every $x\in I_\beta\cap\Q(\beta)$ is eventually periodic. In particular, the expansion of the left boundary point $-\frac{\beta}{\beta+1}$ of the interval $I_\beta$ is eventually periodic. Numbers $\beta$ such that $d_{-\beta}\big(\tfrac{-\beta}{\beta+1}\big)$ is eventually periodic are called Ito-Sadahiro numbers in~\cite{MaPeISnumbers} and Yrrap numbers in~\cite{LiaoSteiner}, where the authors show that the set of Yrrap numbers does not coincide with the set of Parry numbers (see~\cite{Parry}).

Note that the proof of our result, stated as Theorem~\ref{t:obec}, follows the considerations of Akiyama~\cite{akiyamaGreedy} used for positive base systems. The question of admissibility is however less simple for negative bases which is reflected by an additional assumption in the following lemma, justified in Example~\ref{ex:vyjimka}.

\begin{lem}\label{l:obecadmis}
Let $d_{-\beta}\big(\tfrac{-\beta}{\beta+1}\big)$ be eventually periodic not ending in
$0^\omega$ and let $t$ be the length of the longest string of $0$'s
in $d_{-\beta}\big(\tfrac{-\beta}{\beta+1}\big)$. Then for every $(-\beta)$-admissible string
$0a_n\cdots a_00^\omega$, the string $(a_n\cdots
a_00^{t+k})^\omega$ is also admissible for every $k\geq 1$.
\end{lem}

\pfz
Let us verify the conditions of admissibility,
\begin{equation}\label{eq:dk}
d_{-\beta}\big(\tfrac{-\beta}{\beta+1}\big) \preceq_{\text{\tiny alt}} a_i\cdots a_0 (0^{t+k}a_n\cdots
a_0)^\omega \prec_{\text{\tiny alt}} \lim_{\varepsilon\to0+}d_{-\beta}\big(\tfrac{1}{\beta+1}-\varepsilon\big)\,,
\end{equation}
for all $0\leq i\leq n$. Assume for contradiction, that the left
inequality in~\eqref{eq:dk} is not satisfied for some $i$. Then
$$
a_i\cdots a_0 (0^{t+k}a_n\cdots a_0)^\omega \prec_{\text{\tiny
alt}} d_{-\beta}\big(\tfrac{-\beta}{\beta+1}\big) \preceq_{\text{\tiny alt}} a_i\cdots a_0 0^\omega\,.
$$
The strings on the left and right side of the inequality have a
common prefix $a_i\cdots a_0 0^{t+k}$. Necessarily, the string
$d_{-\beta}\big(\tfrac{-\beta}{\beta+1}\big)$ has also such prefix. Therefore $d_{-\beta}\big(\tfrac{-\beta}{\beta+1}\big)$ contains the string
$0^{t+k}$ of the length $t+k>t$, which is a contradiction with the
definition for $t$.

Assume now that the right inequality in~\eqref{eq:dk} is not satisfied for some $i$, i.e.
$$
a_i\cdots a_00^\omega \prec_{\text{\tiny alt}} \lim_{\varepsilon\to0+}d_{-\beta}\big(\tfrac{1}{\beta+1}-\varepsilon\big)
\preceq_{\text{\tiny alt}} a_i\cdots a_0 (0^{t+k}a_n\cdots
a_0)^\omega \,.
$$
By the same considerations, we obtain that $\lim_{\varepsilon\to0+}d_{-\beta}\big(\tfrac{1}{\beta+1}-\varepsilon\big)$ contains a
string of 0 of the length $t+k$. The same is then true for $d_{-\beta}\big(\tfrac{-\beta}{\beta+1}\big)$,
because the two strings appearing in the admissibility condition~\eqref{eq:admis} are in strong connection (cf.~\cite{ItoSadahiro}). This is a contradiction completing the proof. \pfk

Note that the restricting assumption in Lemma~\ref{l:obecadmis} is necessary, as is seen from the following example.

\begin{ex}\label{ex:vyjimka}
Let $d_{-\beta}\big(-\frac{\beta}{\beta+1}\big)=d_1d_2\cdots d_p0^\omega$. The admissibility condition~\eqref{eq:admis} for the digit string $x_1x_2\cdots$ now reads
$$
d_1d_2\cdots d_p0^\omega \preceq_{\text{\tiny alt}} x_ix_{i+1}x_{i+2}\cdots \prec_{\text{\tiny alt}}0d_1d_2\cdots d_p0^\omega\,.
$$
It can be easily seen that the string $1d_1d_2\cdots d_p0^\omega$ is the $(-\beta)$-expansion of $\frac{-1}{\beta(\beta+1)}$, i.e. $01d_1d_2\cdots d_p0^\omega$ is admissible. However,
$(1d_1d_2\cdots d_p0^k)^\omega$ is not admissible for infinitely many values of $k$. In particular,
$$
d_1d_2\cdots d_p0^k (1d_1d_2\cdots d_p0^k)^\omega \prec_{\text{\tiny alt}} d_1d_2\cdots d_p0^\omega
$$
for every odd $k$ if $p$ is odd, and for every even $k$ if $p$ is even.
\end{ex}

\begin{lem}\label{l:obecodhad}
Let $\beta$ be an algebraic unit. Then there exists a positive
constant $c>0$ such that for every
$z\in\Z_{-\beta}\setminus(-\beta)\Z_{-\beta}$ we have
$|z^{(j)}|\geq c$ for at least one of the Galois images
$z^{(2)},\dots,z^{(d)}$ of $z$.
\end{lem}

\pfz
First we show an auxiliary statement:

Let $c_2,\dots,c_d$ be arbitrary positive constants. Then there
exist an integer $p$ such that for every $y\in{\rm Fin}(-\beta)$
with $(-\beta)$-expansion $d_{-\beta}(y)=d_1d_2\cdots d_p0^\omega$,
$d_p\neq 0$, we have $|y^{(j)}|>c_j$ for at least one of the
Galois images $y^{(2)},\dots,y^{(d)}$ of $y$.

We show the statement by contradiction. Suppose that for every
integer $p$ there exists a number ${}^py\in{\rm Fin}(-\beta)$ with
$(-\beta)$-expansion of length $p$ such that for all
$j\in\{2,\dots,d\}$ we have $|{}^py^{(j)}|\leq c_j$. (Note that
the condition on the length of the $(-\beta)$-expansion of ${}^py$
ensures that ${}^py$, $p\in\N$, are distinct.) Since $\beta$ is a
unit, we have ${\rm Fin}(-\beta)\subset \Z[\beta]$, and thus we
have infinitely many different elements ${}^py\in\Z[\beta]$ such
that $|{}^py|\leq 1$, $|{}^py^{(j)}|\leq c_j$. This is impossible
for the following reason: There is a bijection between the set
$\Z[\beta]$ and a $d$-dimensional lattice $\Z^d$, namely
$$
y=y_0+y_1(-\beta)+\cdots y_{d-1}(-\beta)^{d-1}\in\Z[\beta] \mapsto
(y_0,y_1,\dots,y_{d-1})\in\Z^d\,.
$$
Conditions $|{}^py|\leq 1$, $|{}^py^{(j)}|\leq c_j$ are satisfied
by lattice points in a bounded parallelepiped in $\R^d$, and such
lattice points can be only finitely many.

In order to prove the statement of Lemma~\ref{l:obecodhad}, denote
$$
H_j:= \sup\{|z^{(j)}| \mid z\in\Z_{-\beta}\} \leq
\sum_{k=0}^\infty \lfloor\beta\rfloor |\beta^{(j)}|^k =
\frac{\lfloor\beta\rfloor }{1-|\beta^{(j)}|}\,.
$$
We apply the auxiliary statement above with $c_j=2H_j$ for
$j=2,3,\dots,d$ to find the integer $p\geq 1$. Put
$$
c:=\min\Big\{|\beta^{(j)}|^p H_j \,\Big|\, j=2,3,\dots,d\Big\}\,.
$$
Let $z\in\Z_{-\beta}\setminus(-\beta)\Z_{-\beta}$. Then
$\tfrac{z}{(-\beta)^{p}}$ can be written
$$
\frac{z}{(-\beta)^{p}} = w+y\quad\text{ for some $w\in\Z_{-\beta}$
and $y$ with $d_{-\beta}(y)=d_1d_2\cdots d_p0^\omega$, $d_p\neq 0$.}
$$
From the auxiliary fact, there is a $j\in\{2,3,\dots,d\}$ such
that
$$
\bigg|\Big(\frac{z}{(-\beta)^p}\Big)^{(j)}\bigg| \geq |y^{(j)}| -
|w^{(j)}| > c_j - H_j = H_j\,.
$$
Consequently, we have $|z^{(j)}|> |-\beta^{(j)}|^p H_j \geq c$,
which concludes the proof. \pfk

\begin{thm}\label{t:obec}
Let $d_{-\beta}\big(\tfrac{-\beta}{\beta+1}\big)$ be eventually periodic not ending in
$0^\omega$. Let $\beta$ satisfy ${\rm Fin(-\beta)}=\Z[\beta]$.
Then there exists a constant $\gamma>0$ such that any rational
$x\in(-\gamma,\gamma)$ has purely periodic $(-\beta)$-expansion.
\end{thm}

\pfz Let $t$ be the constant from Lemma~\ref{l:obecadmis}, let $c$
be the constant from Lemma~\ref{l:obecodhad} and let $q$ be an
integer such that $\frac{\beta+1}{\beta^q}<1$. Put $\gamma:=
\min\{c,\frac1{\beta^{t+q+1}}\}$.

Consider $x\in(-\gamma,\gamma)\cap\Q$. According to
Lemma~\ref{l:apendix}, there exists an even $N\in\N$ such that
$$
z:=x\big((-\beta)^N-1\big)=x(\beta^N-1)\in\Z[\beta]\,.
$$
Since by assumption, $\beta$ satisfies $\Z[\beta]={\rm
Fin}(-\beta)$, we can write
$$
z=\sum_{i=s}^{n} a_j(-\beta)^j\,,\quad \text{ where } a_n\neq 0,\
a_s\neq 0\,,
$$
and $0a_n\cdots a_s0^\omega$ is an admissible string. We show two
properties of the indices, namely
$$
{\rm (i)}\quad n+t+1<N\qquad\text{ and }\qquad {\rm (ii)}\quad
s\geq 0\,.
$$

For (i), consider first $n$ even, i.e. $n=2k$. Then by
Lemma~\ref{l:intervaly},
$$
\frac{\beta^{n}}{\beta+1} \leq z = x(\beta^N-1) < \gamma\beta^N <
\beta^{N-t-q-1}\,,
$$
which implies
$$
\beta^{n+t+1} < \beta^N \frac{\beta+1}{\beta^q}\leq \beta^N\,,
$$
i.e. $n+t+1<N$. If $n$ is odd, i.e. $n=2k-1$, then again by
Lemma~\ref{l:intervaly}, we have
$$
\frac{\beta^{n}}{\beta+1} \leq -z = -x(\beta^N-1) < \gamma\beta^N
$$
and the same argument implies $n+t+1<N$.

For (ii), put $w=\beta^{-s}z$. Then
$w\in\Z_{-\beta}\setminus(-\beta)\Z_{-\beta}$ and by
Lemma~\ref{l:obecodhad}, there exists a $j\in\{2,3,\dots,d\}$ such
that $|w^{(j)}|\geq c$. We thus have
$$
\gamma\leq c\leq |w^{(j)}| = |\beta^{(j)}|^{-s}|x|
\big(1-(\beta^{(j)})^N\big) < |\beta^{(j)}|^{-s} \gamma\,.
$$
As $\beta$ is a Pisot or Salem number, we have $s\geq 0$. If
$s>0$, we set $a_{s-1}=a_{s-2}=\cdots = a_0=0$.

As $0a_n\cdots a_00^\omega$ is admissible and $N-n-1\geq t+1$, by
Lemma~\ref{l:obecadmis}, the string $(0^{N-n-1}a_n\cdots
a_0)^\omega$ is also admissible. Item (iv) of
Section~\ref{sec:mainidea} then implies that
$x=\frac{z}{(-\beta)^N-1}$ has a purely periodic
$(-\beta)$-expansion, $d_{-\beta}(x)=(0^{N-n-1}a_n\cdots a_0)^\omega$.
 \pfk

\section{Comments}

In this paper we show that under certain conditions, there exists an interval $J\subset I_\beta$ containing 0, such that the $(-\beta)$-expansion of every rational number in $J$ is
purely periodic. Moreover, we describe this interval exactly in case of quadratic bases $\beta$.
Such precise description is possible based on the precise knowledge of the set $\{z'\mid
z\in\Z_{-\beta}\}$. Note that in the case of positive base systems, the parameter $\gamma(\beta)$ is known explicitly for quadratic Pisot numbers $\beta$ and one more example, the minimal Pisot number~\cite{AFSS}.

Note that Theorem~\ref{t:obec} requires besides the finiteness property ${\rm Fin(-\beta)}=\Z[\beta]$
a further assumption that the base $-\beta$ has an eventually periodic but not a finite
expansion $d_{-\beta}\big(\frac{-\beta}{\beta+1}\big)$. It is likely that this assumption does not
represent a real limitation. We have reasons to conjecture that
$\beta$ with finite $d_{-\beta}\big(\frac{-\beta}{\beta+1}\big)=d_1\cdots d_p0^\omega$ never satisfies
${\rm Fin(-\beta)}=\Z[\beta]$, which is also required by the
Theorem. To support such a conjecture, consider the following
example.

\begin{ex}
Consider non-negative integers $d_1,\dots,d_p$ satisfying
$d_1>\max\{d_2,\dots,d_p\}+1$. In~\cite{steiner} one can find a sufficient condition
on a digit string to be the $(-\beta)$-expansion of the left boundary point $l_\beta=\frac{-\beta}{\beta+1}$ for some $\beta$. The string $d_1\dots d_p0^\omega$ automatically satisfies this condition, i.e. there exists $\beta>1$ such that $d_{-\beta}(l_\beta)=d_1\cdots d_p0^\omega$.

One can verify that the
$(-\beta)$-expansion of $\beta-1-d_1\in\Z[\beta]$ is infinite,
$$
d_{-\beta}(\beta-1-d_1) = (d_2+1)(d_3+1)\cdots (d_p+1)1^{\omega}\,.
$$
\end{ex}

As said above, we have shown on the family of quadratic negative bases $-\beta$ with $\beta^2=m\beta+1$ that finiteness property is not necessary for existence of an interval $J$ whose all rational numbers have purely periodic $(-\beta)$-expansion. For such bases $\beta$, we have ${\rm Fin}(-\beta)-{\rm Fin}(-\beta)\not\subset{\rm Fin}(-\beta)$ but we still have
${\rm Fin}(-\beta)+{\rm Fin}(-\beta)\subset{\rm Fin}(-\beta)$.

Recall that for numeration systems with positive bases a necessary and sufficient condition so that $\gamma(\beta)>0$ is known also for $\beta$ cubic. This result uses the description of cubic Pisot numbers with finiteness property~\cite{akiyamaCubic}. No such description is known for negative cubic base. Any result about
closedness of ${\rm Fin}(-\beta)$ under addition or both addition and subtraction could be of help.


\subsubsection*{Acknowledgments.}

We acknowledge financial support by the Czech Science Foundation
grant 201/09/0584 and by the grants MSM6840770039 and LC06002 of
the Ministry of Education, Youth, and Sports of the Czech
Republic.


\end{document}